\documentclass[12pt,a4paper,leqno]{article}

\usepackage{color}
\usepackage[T1]{fontenc}
\usepackage[latin1]{inputenc}
\usepackage[french]{babel}
\usepackage{amsmath}
\usepackage{amssymb}
\usepackage{url}
\usepackage{listings}

\numberwithin{equation}{section}
\newcommand{\suppress}[1]{}

\def\iy{\infty}

\def\al{\alpha}
\def\be{\beta}
\def\ga{\gamma}

\def\de{\delta}

\def\ka{\kappa}

\def\ze{\zeta}

\def\vfi{\varphi}
\def\om{\omega}
\def\Om{\Omega}

\def\ni{\noindent}

\def\LR{\Longrightarrow}

\def\leq{\leqslant}
\def\geq{\geqslant}
\def\mk{\medskip}

\newtheorem{theorem}{Theorem}[section]

\newtheorem{prop}{Proposition}[section]
\newtheorem{coro}{Corollary}[section]
\newtheorem{lem}{Lemma}[section]

\newenvironment{dem}{\smallskip\noindent{\bf Proof}~:%
{\nopagebreak[0]}}%
{\nopagebreak[0]\hfill$\Box$ \medskip}%

\title{Small values of the Euler function and the Riemann hypothesis}
\author{Jean-Louis Nicolas\footnote{Research partially 
supported by CNRS, Institut Camille Jordan, UMR 5208.
This article will be published in Acta Arithmetica.}}

\begin{document}
%\setcounter{page}{0}
%\date{}
\maketitle

\hfill
\begin{minipage}[t]{50mm}
\`A Andr\'e  Schinzel pour son 75\`eme 
anni\-versaire, en tr\`es amical hommage.
\end{minipage}
\medskip

\def\abstractname{Abstract}
\begin{abstract}
Let $\vfi$ be Euler's function, $\ga$ be Euler's constant
and $N_k$ be the product of the first $k$ primes.
In this article, we consider the function 
$c(n) =(n/\vfi(n)-e^\ga\log\log n)\sqrt{\log n}$. Under 
Riemann's hypothesis, it is proved that $c(N_k)$ is bounded and 
explicit bounds are given while, if Riemann's hypothesis fails, 
$c(N_k)$ is not bounded above or below.
\end{abstract}

\medskip

\ni

\medskip

\ni

%%%%%%%%%%%%%%%%%%%%%%%%%%%%%%%%%%%%%%%%%%%%%%%%%%%%%%%%%%
%%%%%%%%%%%%%%%%%%%%%%%%%%%%%%%%%%%%%%%%%%%%%%%%%%%%%%%%%
\section{Introduction}\label{parint}
%%%%%%%%%%%%%%%%%%%%%%%%%%%%%%%%%%%%%%%%%%%%%%%%%%%%%%%%%%
%%%%%%%%%%%%%%%%%%%%%%%%%%%%%%%%%%%%%%%%%%%%%%%%%%%%%%%%%%

Let $\vfi$ be the Euler function. In 1903, it was proved by E. Landau (cf.
\cite[\S 59]{Lan} and \cite[Theorem 328]{HW}) that
\begin{equation*}
\limsup_{n\to \iy} \frac{n}{\vfi(n) \log \log n}=e^\ga=
1.7810724179\ldots  
\end{equation*}
where $\ga=0.5772156649\ldots$     
is Euler's constant.

In 1962, J. B. Rosser and L. Schoenfeld proved (cf. \cite[Theorem 15]{RS})
\begin{equation}\label{RS2.51}
\frac{n}{\vfi(n)} \leq e^\ga\log \log n+\frac{2.51}{\log \log n}
\end{equation}
for $n\geq 3$ and asked if there exists an infinite number of
$n$ such that $n/\vfi(n) > e^\ga\log\log n$. In
\cite{NicJNT}, (cf. also \cite{NicDPP}), I answer this question 
in the affirmative. Soon after, 
A. Schinzel told me that he had worked 
unsuccessfully on this question,  which made me very proud to have
solved it.

For $k\geq 1$, $p_k$ denotes the $k$-th prime and 
$$N_k=2\cdot 3 \cdot 5 \; \ldots\; p_k$$
the primorial number of order $k$. In \cite{NicJNT}, it is proved that
the Riemann hypothesis (for short RH) is equivalent to 
$$\forall k \geq 1, \qquad \frac{N_k}{\vfi(N_k)} >
e^\ga \log \log N_k.$$
The aim of the present paper is to make more precise the results of
\cite{NicJNT} by estimating the quantity
\begin{equation}\label{c}
c(n)=\left(\frac{n}{\vfi(n)} -e^\ga \log \log n \right)\sqrt{\log n}.
\end{equation}
Let us denote by $\rho$ a generic root of the Riemann $\ze$ function
satisfying $0< \Re \rho < 1$. Under RH, $1-\rho=\overline{\rho}$.
It is convenient to define (cf. \cite[p. 159]{Edw})
\begin{equation}\label{beta}
\be=\sum_\rho \frac{1}{\rho(1-\rho)}=
2+\ga-\log\pi-2\log 2=0.0461914179\ldots   
\end{equation}
We shall prove
\begin{theorem}\label{thmRH}
Under the Riemann hypothesis (RH) we have
\begin{equation}\label{rhlimsup}
\limsup_{n\to \iy} c(n)=e^\ga(2+\be)=3.6444150964\ldots     
\end{equation}
\begin{equation}\label{rh120}
\forall n \geq N_{120569}=2\cdot 3 \cdot\ldots\cdot 1591883, 
\qquad c(n) <  e^\ga(2+\be).
\end{equation}
\begin{equation}\label{rh66}
\forall n \geq 2, \qquad c(n) \leq
c(N_{66})=c(2\cdot 3 \cdot\ldots\cdot 317)=
4.0628356921\ldots    
\end{equation}
\begin{equation}\label{rh2}
\forall k \geq 1, \qquad c(N_k) \geq
c(N_{1})=c(2)=2.2085892614\ldots   
\end{equation}
\end{theorem}
We keep the notation of \cite{NicJNT}. For a real $x\geq 2$,  the
usual Chebichev's functions are denoted by
\begin{equation}\label{Che}
\theta(x)=\sum_{p\leq x} \log p\qquad \text{ and } \qquad
\psi(x)=\sum_{p^m\leq x} \log p.
\end{equation}
We set
\begin{equation}\label{f}
f(x)=e^\ga \log \theta(x) \prod_{p\leq x} (1-1/p).
\end{equation}
Mertens's formula yields $\lim_{x\to\iy} f(x)=1$.
In \cite[Th. 3 (c)]{NicJNT} it is shown that, if RH fails, there
exists $b$, $0 < b < 1/2$, such that 
\begin{equation}\label{Omeg}
\log f(x) = \Om_{\pm}(x^{-b}).
\end{equation}
For $p_k \leq x < p_{k+1}$, we have
$f(x)=e^\ga \log \log(N_k) \frac{\vfi(N_k)}{N_k}\cdot$ 
When $k\to\iy$, by observing that the Taylor development
in neighborhood of $1$ yields $\log f(p_k)\sim f(p_k)-1$, we get
$$\log f(p_k) \sim f(p_k)-1=\frac{\vfi(N_k)}{N_k}
\frac{c(N_k)}{\sqrt{\log N_k}} \sim
\frac{e^{-\ga}}{\log \log N_k}\frac{c(N_k)}{\sqrt{\log N_k}},$$
 and it follows from \eqref{Omeg} that, if RH does not hold, then
$$\liminf_{n\to\iy} \;c(n)=-\iy\quad \text{ and } 
\quad \limsup_{n\to\iy}\; c(n)=+\iy.$$
Therefore, from Theorem \ref{thmRH}, we deduce~:
\begin{coro}\label{coroRH}
Each of the four assertions  \eqref{rhlimsup}, \eqref{rh120}, 
\eqref{rh66}, \eqref{rh2} is equivalent to the Riemann hypothesis.
\end{coro}

%%%%%%%%%%%%%%%%%%%%%%%%%%%%%%%%%%%%%%%%%%%%%%%%%%%%%%%%%%
\subsection{Notation and  results used}\label{parNota}
%%%%%%%%%%%%%%%%%%%%%%%%%%%%%%%%%%%%%%%%%%%%%%%%%%%%%%%%%%

If $\theta(x)$ and $\psi(x)$ are the Chebichev functions 
defined by \eqref{Che}, we set
\begin{equation}\label{Sx}
R(x)=\psi(x)-x \qquad \text{ and } \qquad S(x)=\theta(x)-x.
\end{equation}
Under RH, we shall use the upper bound (cf. \cite[(6.3)]{Sch76})
\begin{equation}\label{Tx}
x \geq 599 \quad \LR \quad |S(x)| \leq 
T(x)\stackrel{def}{=\!=} \frac{1}{8\pi}\sqrt x\log ^2 x
\end{equation}
P. Dusart (cf. \cite[Table 6.6]{Dus}) has shown that
\begin{equation}\label{dus}
\theta(x) < x \text{ for } x \le 8\cdot10^{11}
\end{equation}
thus improving the result of R. P. Brent who has
checked \eqref{dus} for $x < 10^{11}$ 
(cf. \cite[p. 360]{Sch76}).
We shall also use (cf. \cite[Theorem 10]{RS}
\begin{equation}\label{th4/5}
\theta(x) \geq 0.84\;x\geq \frac 45 x\qquad \text{ for } x\geq 101.
\end{equation}

\mk

As in \cite{NicJNT}, we define the following integrals
\begin{equation}\label{K}
K(x)=\int_x^\iy \frac{S(t)}{t^2}\left( \frac{1}{\log
    t}+\frac{1}{\log^2 t}\right) dt,
\end{equation}
\begin{equation}\label{J}
J(x)=\int_x^\iy \frac{R(t)}{t^2}\left( \frac{1}{\log
    t}+\frac{1}{\log^2 t}\right) dt,
\end{equation}
and, for $\Re (z) < 1$,
\begin{equation}\label{Fz}
F_z(x)=\int_x^\iy t^{z-2}\left( \frac{1}{\log
    t}+\frac{1}{\log^2 t}\right) dt.
\end{equation}
We also set for $x\geq 1$
\begin{equation}\label{W}
W(x)=\sum_\rho \frac{x^{i\,\Im(\rho)}}{\rho(1-\rho)}
\end{equation}
so that, under RH, from \eqref{beta} we have
\begin{equation}\label{Wbeta}
|W(x)|\leq \be =\sum_\rho \frac{1}{\rho(1-\rho)}\cdot
\end{equation}
We often implicitly use the following result~: for $a$ and $b$
positive, the function
\begin{equation}\label{fab}
t\mapsto \frac{\log^a t}{t^b} \quad \text{is decreasing for}\quad t > e^{a/b}
\end{equation}
and
\begin{equation}\label{fabmax}
\max_{t \geq 1} \frac{\log^a t}{t^b} =\left(\frac{a}{e\, b}\right)^a.
\end{equation}

%%%%%%%%%%%%%%%%%%%%%%%%%%%%%%%%%%%%%%%%%%%%%%%%%%%%%%%%%%
\subsection{Organization of the article}\label{parOrg}
%%%%%%%%%%%%%%%%%%%%%%%%%%%%%%%%%%%%%%%%%%%%%%%%%%%%%%%%%%

In Section \ref{parfx}, the results of \cite{NicJNT} about $f(x)$ 
are revised so as to get effective upper and lower bounds for both
$\log f(x)$ and $1/f(x)-1$ under RH (cf.  Proposition \ref{propf}).

In Section \ref{parcn}, we study $c(N_k)$ and $c(n)$ in terms of $f(p_k)$.

Section \ref{parproof} is devoted to the proof of Theorem  \ref{thmRH}.

%%%%%%%%%%%%%%%%%%%%%%%%%%%%%%%%%%%%%%%%%%%%%%
%%%%%%%%%%%%%%%%%%%%%%%%%%%%%%%%%%%%%%%%%%%%%%%%%%%%%%%%%
\section{Estimate of $\log(f(x))$}\label{parfx}
%%%%%%%%%%%%%%%%%%%%%%%%%%%%%%%%%%%%%%%%%%%%%%%%%%%%%%%%%%
%%%%%%%%%%%%%%%%%%%%%%%%%%%%%%%%%%%%%%%%%%%%%%%%%%%%%%%%%%%
The following lemma is Proposition 1 of \cite{NicJNT}.
\begin{lem}\label{lemKx}
For $x\geq 121$, we have
\begin{equation}\label{Kineg}
K(x)-\frac{S^2(x)}{x^2\log x} \leq \log f(x) \leq K(x)+\frac{1}{2(x-1)}.
\end{equation} 
\end{lem}
The next lemma is a slight improvement of Lemma 1 of \cite{NicJNT}.
\begin{lem}\label{lemFz}
Let $x$ be a real number, $x > 1$. For $\Re z < 1$, we have
\begin{equation}\label{rz}
F_z(x)=\frac{x^{z-1}}{(1-z)\log x}+r_z(x)\;\text{ with } \; 
r_z(x)=\int_x^\iy -\frac{z t^{z-2}}{(1-z)\log^2 t} dt
\end{equation} 
and, if $\Re z=1/2$,
\begin{equation}\label{rzabs}
|r_z(x)|\leq\frac{1}{|1-z|\sqrt x \log^2 x} \left(1+\frac{4}{\log x}\right).
\end{equation} 
Moreover, for $z=1/2$, we have
\begin{equation}\label{Fdemi}
\frac{2}{\sqrt x \log x}- \frac{2}{\sqrt x \log^2 x}
\leq F_{1/2}(x)\leq \frac{2}{\sqrt x \log x} 
-\frac{2}{\sqrt x \log^2 x} + \frac{8}{\sqrt x \log^3 x} 
\end{equation} 
and, for $z=1/3$, 
\begin{equation}\label{Ftiers}
0 \leq F_{1/3}(x)\leq \frac{3}{2x^{2/3} \log x} \cdot
\end{equation} 
\end{lem}

\begin{dem}
The proof of \eqref{rz} is easy by taking the derivative. By 
partial summation,
we get
\begin{equation}\label{rpart}
r_z(x)=-\frac{z}{1-z} \left(\frac{x^{z-1}}{(1-z)\log^2 x}
+\int_x^\iy\frac{2\ t^{z-2}}{(z-1)\log^3 t} dt\right).
\end{equation} 
If we assume  $\Re z=1/2$, we have $1-z=\overline{z}$ and
$$|r_z(x)|\leq\frac{1}{|1-z|\sqrt x \log^2 x} +\frac{2}{|1-z|\log^3 x}
\int_x^\iy t^{-3/2} dt$$
which yields \eqref{rzabs}.
The proof of \eqref{Fdemi} follows from \eqref{rz} and \eqref{rpart}
by choosing $z=1/2$. The proof of \eqref{Ftiers} follows from
\eqref{rz} since $r_{1/3}$ is negative.
\end{dem}

To estimate the difference $J(x)-K(x)$, we need  Lemma \ref{lempsith} which,
under RH, is an improvement of Propositions 3.1 and 3.2 of \cite{Dus} 
(obtained without assuming RH). The following lemma will be useful for
proving Lemma \ref{lempsith}.

\begin{lem}\label{lemHL}
Let $\ka=\ka(x)=\lfloor\frac{\log x}{\log 2}\rfloor$ the largest
integer such that $x^{1/\ka}\geq 2$. For $x\geq 16$, we set
$$H(x) =1+\sum_{k=4}^\ka x^{1/k-1/3}$$
and for $x\geq 4$
$$L(x) = \sum_{k=2}^\ka \ell_k(x) \quad \text{ with }\quad 
\ell_k(x)=\frac{T(x^{1/k})}{x^{1/3} }=
\frac{\log^2 x }{8\pi\,k^2  x^{1/3-1/(2k)}}\cdot$$
(i) For $j \geq 9$ and $x\geq 2^j$, 
$H(x) \leq H(2^j)$ holds.

\ni
(ii) For $j\geq 35$ and $x\geq 2^j$, 
$L(x) \leq L(2^j)$ holds.
\end{lem}

\begin{dem}
The function $H$ is continuous and decreasing on $[2^j,2^{j+1})$; so, 
to show (i), it suffices to prove for $j\geq 9$
\begin{equation}\label{H2j}
H(2^j) \geq H(2^{j+1}).
\end{equation} 
If $9\leq j\leq 19$, we check \eqref{H2j} by computation. If $j\geq
20$, we have
\begin{eqnarray*}
H(2^j) -H(2^{j+1}) &=& \sum_{k=4}^j 2^{j\left(\frac 1k -\frac 13\right)} 
\left(1-2^{\left(\frac 1k -\frac 13\right)} \right) -
2^{(j+1)\left(\frac{1}{j+1} -\frac 13\right)} \\
&\geq& 2^{j\left(\frac 14-\frac 13\right)} 
\left(1-2^{\left(\frac 14 -\frac 13\right)} \right) -
2^{(j+1)\left(\frac{1}{j+1} -\frac 13\right)} \\
&=& 2^{-\frac j3}\left[(1-2^{-\frac{1}{12}})2^\frac j4-2^\frac 23\right]
\end{eqnarray*}
which proves \eqref{H2j} since the above bracket is $\geq
(1-2^{-\frac{1}{12}})2^\frac{20}{4}-2^\frac 23=0.208\ldots$
and therefore positive.

Let us assume that $j\geq 35$ so that $2^j\geq e^{24}$ holds. From
\eqref{fab}, for each $k\geq 2$, $x\mapsto \ell_k(x)$ is decreasing
for $x\geq 2^j$ so that $L$ is decreasing on $[2^j,2^{j+1})$ and, to show (ii), it
suffices to prove 
\begin{equation}\label{L2j}
L(2^j) \geq L(2^{j+1}).
\end{equation} 
We have
\begin{eqnarray*}
L(2^j) -L(2^{j+1}) &=& \sum_{k=2}^j
\left\{\ell_k(2^j)-\ell_k(2^{j+1})\right\}    -\ell_{j+1}(2^{j+1})\\
&\geq& \ell_2(2^j)-\ell_2(2^{j+1})-\ell_{j+1}(2^{j+1})\\
&=& \frac{\log^2 2}{32\pi}2^{-\frac j3}\left\{ 2^\frac
  j4\left[j^2-2^{-\frac{1}{12}}(j+1)^2\right]-4\cdot 2^\frac 16\right\}.
\end{eqnarray*}
For $j\geq \frac{1}{2^{1/12}-1}=16.81\ldots$,
the above square bracket is increasing on $j$ 
and it is positive for $j=35$. Therefore, the curly
bracket is increasing for $j \geq 35$ and, since its 
value for $j=35$ is equal to $744.17\ldots$,
\eqref{L2j} is proved for $j\geq 35$. 
\end{dem}

\begin{lem}\label{lempsith}
Under RH, we have
\begin{equation}\label{psithmin}
\psi(x)-\theta(x) \geq \sqrt x, \qquad \text{ for } x \geq 121
\end{equation} 
and, for $x \geq 1$,
\begin{equation}\label{psithmax}
\frac{\psi(x)-\theta(x) -\sqrt x}{x^{1/3}}\leq 1.332768\ldots \leq \frac 43\cdot
\end{equation} 
\end{lem}

\begin{dem}
For $x < 599^3$, we check \eqref{psithmin} by computation.
Note that $599$ is prime. Let
$q_0=1$, and let $q_1=4,q_2=8,q_3=9,\ldots, q_{1922}=599^3$ 
be the sequence of powers (with exponent
$\geq 2$) of primes not exceeding $599^3$. On the intervals
$[q_i,q_{i+1})$, the function $\psi-\theta$ is constant and
$x\mapsto (\psi(x)-\theta(x))/\sqrt x$ is decreasing. For $11\leq i \leq
1921$ (i.e. $121 \leq q_i < q_{i+1} \leq 599^3$), we calculate
$\de_i=(\psi(q_i)-\theta(q_i))/\sqrt{q_{i+1}}$ and find that 
$\min_{11\leq i \leq 1921}\de_i =\de_{1886}=1.0379\ldots$ 
($q_{1886}=206468161=14369^2$) while $\de_{10} =0.9379\ldots < 1$ 
($q_{10}=81$). 

Now, we assume $x \geq 599^3$, so that, by \eqref{Tx}, we have
\begin{equation}\label{psithmin1}
\psi(x)-\theta(x) \geq \theta(x^{1/2})+\theta(x^{1/3})\geq 
x^{1/2} +x^{1/3}-T(x^{1/2}) -T(x^{1/3}).
\end{equation} 
By using \eqref{fabmax}, we get
$$\frac{T(x^{1/2})}{x^{1/3}}+\frac{T(x^{1/3})}{x^{1/3}}=
\frac{1}{8\pi}\left(\frac{\log^2 x}{4x^{1/12}}+\frac{\log^2
    x}{9x^{1/6}}\right)\leq \frac{20}{\pi e^2}=0.86157\ldots$$
which, with \eqref{psithmin1}, implies
\begin{equation}\label{psithmin2}
\psi(x)-\theta(x) \geq \sqrt x +\left(1-\frac{20}{\pi e^2} \right)x^{1/3}
\geq \sqrt x.
\end{equation} 

The inequality \eqref{psithmax} is Lemma 3 of \cite{robin}. We give
below another proof by considering three cases
according to the values of $x$.

\mk
\ni
{\bf Case 1, $1\leq x < 2^{32}$.} The largest $q_i$ smaller than
$2^{32}$ is $q_{6947}=4293001441=65521^2$. 
On the intervals $[q_i,q_{i+1})$, the function 
$$G(x) \;\stackrel{def}{=\!=} \;\frac{\psi(x)-\theta(x)-\sqrt x}{x^{1/3}}$$
is decreasing. By computing $G(q_0),G(q_1),\ldots,G(q_{6947})$ we get
$$G(x)\leq G(q_{103})=1.332768\ldots\qquad
\qquad [q_{103}=80089=283^2].$$   

\mk
\ni
{\bf Case 2, $2^{32} \leq x < 64\cdot 10^{22}$.}
By using \eqref{dus}, we get
$$\psi(x)-\theta(x)=\sum_{k=2}^\ka \theta(x^{1/k})\leq \sum_{k=2}^\ka
x^{1/k}$$
so that Lemma \ref{lemHL} implies
$G(x) \leq H(x)\leq H(2^{32})=1.31731\ldots$

\mk
\ni
{\bf Case 3, $x\geq 64\cdot 10^{22}  \geq 2^{79}$.}
By \eqref{Tx} and  \eqref{dus}, we get
$$\psi(x)-\theta(x)=\sum_{k=2}^\ka \theta(x^{1/k})\leq \sum_{k=2}^\ka
\left\{x^{1/k}+T(x^{1/k})\right\},$$
whence, from Lemma \ref{lemHL},
$G(x) \leq H(x)+L(x)\leq H(2^{79})+L(2^{79})=1.32386\ldots$
\end{dem}

\begin{coro}\label{coroJK}
For $x\geq 121$, we have
\begin{equation}\label{JK}
F_{1/2}(x) \leq J(x)-K(x) \leq F_{1/2}(x) +\frac 43 F_{1/3}(x).
\end{equation}
\end{coro}
The following lemma is an improvement of \cite[Proposition 2]{NicJNT}. 
\begin{lem}\label{lemJ}
Let us assume that RH holds. For $x > 1$, we may write
\begin{equation}\label{WJ1J2}
J(x)=-\frac{W(x)}{\sqrt x \log x} -J_1(x)-J_2(x)
\end{equation}
with
\begin{equation}\label{J1J2}
0 < J_1(x) \leq \frac{\log (2\pi)}{x\log x}\quad \text{ and } \quad
|J_2(x)| \leq \frac{\be}{\sqrt x \log ^2 x}\left(1+\frac{4}{\log x}\right).
\end{equation}
\end{lem}

\begin{dem}
In \cite[(17)--(19)]{NicJNT}, for $x > 1$, it is proved that
$$J(x)=-\sum_\rho \frac1\rho F_\rho(x)-J_1(x)$$
with $J_1$ satisfying $0 < J_1(x) \leq
\frac{\log (2\pi)}{x\log x}\cdot$ 

Now, by Lemma \ref{lemFz}, we have
$F_\rho(x)=\frac{x^{\rho-1}}{(1-\rho)\log x}+r_\rho(x)$ which yields
\eqref{WJ1J2} by setting $J_2(x)=\sum_\rho \frac1\rho r_\rho(x)$. 
Further, from \eqref{rzabs} and \eqref{beta}, we get the upper bound
for $|J_2(x)|$ given in \eqref{J1J2}.
\end{dem}

\begin{prop}\label{propf}
Under RH, for $x\geq x_0=10^9$, we have
\begin{equation}\label{logf}
  -\frac{2+W(x)}{\sqrt x \log x} +\frac{0.055}{\sqrt x \log^2 x}
\leq \log f(x) \leq 
-\frac{2+W(x)}{\sqrt x \log x} +\frac{2.062}{\sqrt x \log^2 x}
\end{equation} 
and
\begin{equation}\label{1/fx}
 \frac{2+W(x)}{\sqrt x \log x} -\frac{2.062}{\sqrt x \log^2 x}
\leq \frac{1}{f(x)} -1 \leq 
\frac{2+W(x)}{\sqrt x \log x} -\frac{0.054}{\sqrt x \log^2 x}\cdot
\end{equation} 
\end{prop}

\begin{dem}
By collecting the information from
\eqref{Kineg}, \eqref{Tx}, \eqref{JK}, \eqref{WJ1J2}, \eqref{J1J2}, 
\eqref{Fdemi} and \eqref{Ftiers}, for $x\geq 599$, we
get
\begin{eqnarray}\label{mlogf}
\log f(x) &\geq& -\frac{W(x)+2}{\sqrt x \log x}
+\frac{2-\be}{\sqrt x\log^2 x} -\frac{8+4\be}{\sqrt x \log^3 x}
\notag\\
& & \hspace{20mm}-\frac{\log(2\pi)}{x\log x} -
\frac{2}{x^{2/3}\log x}-\frac{\log^3 x}{64\pi^2 \,x}
\end{eqnarray}
and
\begin{equation}\label{Mlogf}
\log f(x) \leq -\frac{W(x)+2}{\sqrt x \log x}+\frac{2+\be}{\sqrt
  x\log^2 x} +\frac{4\be}{\sqrt x \log^3 x}+\frac{1}{2(x-1)}\cdot
\end{equation} 
Since $x\geq x_0=10^9$ holds, \eqref{mlogf} and \eqref{Mlogf} imply
respectively
\begin{eqnarray}\label{mlogf0}
\log f(x) &\geq& -\frac{W(x)+2}{\sqrt x \log x}
 +\frac{1}{\sqrt x\log^2 x} \left(2-\be-\frac{8+4\be}{\log x_0}\right.\notag\\
& &\hspace{20mm}\left.-\frac{\log(2\pi)\log x_0}{\sqrt{x_0}} -\frac{2\log x_0}{x_0^{1/6}}-\frac{\log^5
  x_0}{64\pi^2 \sqrt{x_0}} \right)
\end{eqnarray}
and
\begin{equation}\label{Mlogf0}
\log f(x) \leq -\frac{W(x)+2}{\sqrt x \log x}+\frac{1}{\sqrt x \log^2 x} \left(
2+\be +\frac{4\be}{\log x_0}+\frac{\sqrt{x_0}\log^2
  x_0}{2(x_0-1)}\right)
\end{equation} 
which prove \eqref{logf}.

Setting $v=-\log f(x)$, it follows from \eqref{logf}, \eqref{Wbeta} and
\eqref{beta} that 
$$v\leq \frac{W(x)+2}{\sqrt x \log x}\leq \frac{2+\be}{\sqrt x \log
  x}\leq v_0\;\stackrel{def}{=\!=} \;\frac{2+\be}{\sqrt x_0 \log x_0}=0.00000312\ldots$$
By Taylor's formula, we have $e^v-1\geq v$ (which, with \eqref{logf},
provides the lower bound of \eqref{1/fx}) and 
$$e^v -1-v \leq \frac{e^{v_0}}{2}v^2 \leq
\frac{e^{v_0}(2+\be)^2}{2x\log^2 x} \leq
\frac{e^{v_0}(2+\be)^2}{2\sqrt{x_0}\sqrt x\log^2 x}
=\frac{0.0000662\ldots}{\sqrt x\log^2 x}$$
(which implies the upper bound in \eqref{1/fx}).
\end{dem}

%%%%%%%%%%%%%%%%%%%%%%%%%%%%%%%%%%%%%%%%%%%%%%
%%%%%%%%%%%%%%%%%%%%%%%%%%%%%%%%%%%%%%%%%%%%%%%%%%%%%%%%%
\section{Bounding $c(n)$}\label{parcn}
%%%%%%%%%%%%%%%%%%%%%%%%%%%%%%%%%%%%%%%%%%%%%%%%%%%%%%%%%%
%%%%%%%%%%%%%%%%%%%%%%%%%%%%%%%%%%%%%%%%%%%%%%%%%%%%%%%%%%%

\begin{lem}\label{lemcn}
Let $n$ and $k$ be two integers satisfying $n\geq 2$ and $k \geq
1$. Let us assume that either the 
number $j=\om(n)$ of distinct prime factors of $n$ is equal to $k$ or
that $N_k \leq n < N_{k+1}$ holds. We have
\begin{equation}\label{cn<}
c(n) \leq c(N_k).
\end{equation}
\end{lem}

\begin{dem}
It follows from our hypothesis that $n\geq N_k$ and $j \leq k$ hold.
Let us write
$n=q_1^{\al_1}q_2^{\al_2}\ldots q_j^{\al_j}$ (with $q_1 < q_2 < \ldots
< q_j$ as defined in the proof of Lemma \ref{lempsith}). We have
$$\frac{n}{\vfi(n)}=\prod_{i=1}^j \frac{1}{1-1/q_i} \leq 
\prod_{i=1}^j \frac{1}{1-1/p_i} \leq \prod_{i=1}^k \frac{1}{1-1/p_i} 
=\frac{N_k}{\vfi(N_k)}$$ 
which yields
\begin{equation}\label{cn<hn}
c(n) \leq \left(\frac{N_k}{\vfi(N_k)}-e^\gamma \log \log
  n\right)\sqrt{\log n} \;\stackrel{def}{=\!=} \;h(n)
\end{equation}
and $h(n)$ can be extended to a real number $n$. Further,
\begin{eqnarray*}
\frac{d}{dn} h(n) &=& \frac{1}{2n\sqrt{\log n}}
\left(\frac{N_k}{\vfi(N_k)}-e^\ga \log \log  n-2e^\ga\right)\\
&\leq& \frac{1}{2n\sqrt{\log n}}\left(\frac{N_k}{\vfi(N_k)}-
e^\ga \log \log  N_k-2e^\ga\right).
\end{eqnarray*}
If $k=1$ or $2$, it is easy to see that the above parenthesis is
negative, while, if $k\geq 3$, by \eqref{RS2.51}, it is smaller than
$\frac{2.51}{\log \log N_k}-2e^\ga$ which is also negative because 
$\log \log N_k \geq \log \log 30=1.22\ldots$ 
Therefore, we get $h(n) \leq
h(N_k)=c(N_k)$, which, with \eqref{cn<hn}, completes the proof of
Lemma \ref{lemcn}.  
\end{dem}

\begin{prop}\label{propcn}
Let us assume that $x_0=10^9 \leq p_k\leq x < p_{k+1}$ holds.
Under RH, we have
\begin{equation}\label{cNk<}
c(N_k) \leq e^\ga(2+W(x))-\frac{0.07}{\log x} \leq 
e^\ga(2+\be)-\frac{0.07}{\log x} 
\end{equation}
and
\begin{equation}\label{cNk>}
c(N_k) \geq e^\ga(2+W(x))-\frac{3.7}{\log x} \geq 
e^\ga(2-\be)-\frac{3.7}{\log x} \cdot
\end{equation}
\end{prop}

\begin{dem}
From \eqref{c} and \eqref{f}, we get
\begin{equation}\label{cN1}
c(N_k)=e^\ga \sqrt{\theta(x)}\log \theta(x)\left(\frac{1}{f(x)}-1\right)\cdot
\end{equation}
By the fundamental theorem of calculus, \eqref{th4/5} and \eqref{Tx}, we have
\begin{eqnarray*}
|\sqrt{\theta(x)}\log \theta(x)-\sqrt x\log x|
&=& \left|\int_x^{\theta(x)}\frac{\log t+2}{2\sqrt t} dt\right| \leq
|\theta(x)-x| \frac{\log(4x/5)+2}{2\sqrt{4x/5}}\\
&\leq& \frac{\sqrt 5}{4}T(x) \frac{\log x+2}{\sqrt x}=
\frac{\sqrt 5}{32\pi}\log^2 x(\log x+2)
\end{eqnarray*}
whence
$$\left|\frac{\sqrt{\theta(x)}\log \theta(x)}{\sqrt x\log x}-1\right|\! \leq\!
\frac{\sqrt 5 \log^2 x(\log x+2)}{32\pi\sqrt x \log x}\!\leq\!
\frac{\sqrt 5 \log^2 x_0(\log x_0+2)}{32\pi\sqrt{x_0} \log x}\!\leq\!
\frac{0.0069}{\log x}\cdot$$
Therefore, \eqref{cN1}, \eqref{1/fx} and \eqref{Wbeta} yield
\begin{eqnarray*}
c(N_k) &\leq& e^\ga\left(2+W(x)-\frac{0.054}{\log x}\right)
\left(1+\frac{0.0069}{\log x} \right)\\
&\leq& e^\ga(2+W(x))-\frac{e^\ga}{\log x}(0.054-0.0069(2+\be))
\end{eqnarray*}
which proves \eqref{cNk<}. The proof of \eqref{cNk>} is similar. 
\end{dem}

%%%%%%%%%%%%%%%%%%%%%%%%%%%%%%%%%%%%%%%%%%%%%%
%%%%%%%%%%%%%%%%%%%%%%%%%%%%%%%%%%%%%%%%%%%%%%%%%%%%%%%%%
\section{Proof of Theorem \ref{thmRH}}\label{parproof}
%%%%%%%%%%%%%%%%%%%%%%%%%%%%%%%%%%%%%%%%%%%%%%%%%%%%%%%%%%
%%%%%%%%%%%%%%%%%%%%%%%%%%%%%%%%%%%%%%%%%%%%%%%%%%%%%%%%%%%

It follows from \eqref{cn<}, \eqref{cNk<} and  \eqref{cNk>}  that
$$\limsup_{n\to\iy}c(n)=e^\ga(2+\limsup_{x\to\iy} W(x)).$$
As observed in \cite[p. 383]{NicJNT}, 
by the pigeonhole principle (cf.  \cite[\S 2.11]{EMF} or \cite[\S
11.12]{HW}), one can show that 
$\limsup_{x\to\iy} W(x)=\be$, which proves \eqref{rhlimsup}.

To show the other points of Theorem \ref{thmRH}, we first consider
$k_0=50847534$, the number of primes up to $x_0=10^9$.
For all $k\leq k_0$, we have calculated $c(N_k)$ in Maple with $30$ 
decimal digits, so that we may think that the first ten are correct.

We have found that for $k_1=120568 < k \leq k_0$, $c(N_k) <
e^\ga(2+\be)$ holds (while $c(N_{k_1}) =3.6444180\ldots >  e^\ga(2+\be)$)
and for $1\leq k \leq k_0$, we have
$c(N_1)=c(2)\leq c(N_k) \leq c(N_{66})$.

Further, for $k > k_0$, \eqref{cNk<} implies $c(N_k) < e^\ga(2+\be) <
c(N_{66})$ which, together with Lemma \ref{lemcn}, proves \eqref{rh120}
and \eqref{rh66}.

As a challenge, for $k_1=120568$, I ask to find the largest number $M$
such that $M < N_{k_1+1}$ and $c(M) \geq e^\ga(2+\be)$. Note
that $M > N_{k_1}$ holds since, for $n=N_{k_1-1}p_{k_1+1}$, we have 
$c(n)=3.6444178\ldots > e^\ga(2+\be)$. Another challenge is to
determine all the $n$'s satisfying $n < N_{k_1+1}$ and $c(n) >
e^\ga(2+\be)$.

 Finally, for $k > k_0$, \eqref{cNk>} implies 
$$c(N_k)\geq e^\ga(2-\be)-\frac{3.7}{\log (10^9)}=3.30\ldots > c(2)$$
which completes the proof of \eqref{rh2} and of Theorem
\ref{thmRH}. 
\hfill  $\Box$

\mk
It is not known if $\liminf_{x\to\iy} W(x)=-\be$. Let
$\rho_1=1/2+ i\ t_1$ with
$t_1=14.13472\ldots$ the first zero of $\zeta$. By using
a theorem of Landau (cf. \cite[Th. 6.1 and \S 2.4]{EMF}), it is
possible to prove that $\liminf_{x\to\iy} W(x)\leq
-1/(\rho_1(1-\rho_1))=-0.00499\ldots$ A smaller upper bound is wanted.
 
An interesting question is the following~: assume that RH fails. Is it
possible to get an upper bound for $k$ such that $k > k_0$ and either 
$c(N_k) > e^\ga(2+\be)$ or $c(N_k) < c(2)$?

\bigskip

\ni
{\bf Acknowledgements.} We would like to thank the anonymous referee
for his or her careful reading of our article and his or her valuable suggestions.

\vspace{-5mm}

\def\refname{References}

\ni
Jean-Louis Nicolas,
Université de Lyon, Université de Lyon 1, CNRS,\\
Institut Camille Jordan, Math\'ematiques, Bât. Doyen Jean Braconnier,\\
Université Claude Bernard (Lyon 1), 21 Avenue Claude Bernard,\\
F-69622 Villeurbanne c\'edex, France.

\ni
\url{jlnicola@in2p3.fr}, \quad 
\url{http://math.univ-lyon1.fr/~nicolas/}
\end{document}